\documentclass{amsart}

\usepackage {amssymb}
\usepackage{verbatim}
\usepackage[all]{xy}
\usepackage{graphicx}
\usepackage[usenames]{color}

\newcommand{\Hom}{\mathrm{Hom}}
\renewcommand{\max}{\mathrm{max}}

\newtheorem{theorem}[equation]{Theorem}

\newtheorem{cor}[equation]{Corollary}
\newtheorem{prop}[equation]{Proposition}
\newtheorem{conj}[equation]{Conjecture}

\theoremstyle{definition} 
\newtheorem{defn}[equation]{Definition}
\newtheorem{example}[equation]{Example}
\newtheorem{remark}[equation]{Remark}

\theoremstyle{remark}

\numberwithin{equation}{section}
\numberwithin{figure}{section}

\newcommand{\SL}{\mathrm{SL}}

\newcommand{\A}{\mathbb{A}}

\newcommand{\TQB}{T^*_Q\caB}
\newcommand{\TQP}{T^*_Q\caP}

\renewcommand{\top}{{\mathrm{top}}}
\newcommand{\ol}{\overline}

\newcommand{\field}{\mathbb}
\newcommand{\liealgebra}{\mathfrak}
\newcommand{\la}{\liealgebra}

\newcommand{\C}{{\field C}}
\newcommand{\R}{{\field R}}

\renewcommand{\u}{{\liealgebra u}}

\newcommand{\n}{{\la n}}

\newcommand{\AV}{\mathrm{AV}}
\newcommand{\CV}{\mathrm{CV}}
\newcommand{\av}{\mathrm{av}}
\newcommand{\cv}{\mathrm{cv}}
\newcommand{\Ann}{\mathrm{Ann}}

\newcommand{\gr}{\mathrm{gr}}

\renewcommand{\H}{\mathrm{H}}
\newcommand{\bfM}{\mathbf{M}}

\newcommand{\triv}{1 \! \! 1}
\newcommand{\ind}{\mathrm{ind}}
\newcommand{\sgn}{\mathrm{sgn}}
\newcommand{\lra}{\longrightarrow}

\newcommand{\wt}{\widetilde}

\newcommand{\bs}{\backslash}

\newcommand{\beq}{\begin{equation}}
\newcommand{\eeq}{\end{equation}}

\newcommand{\std}{\textrm{std}}

\newcommand{\Irr}{\mathrm{Irr}}
\newcommand{\muinv}{\mu^{-1}}

\newcommand{\supp}{\mathrm{supp}_{\circ}}
\newcommand{\U}{\mathrm{U}}
\newcommand{\GL}{\mathrm{GL}}
\newcommand{\Sp}{\mathrm{Sp}}
\newcommand{\SO}{\mathrm{SO}}

\newcommand{\frb}{\mathfrak{b}}

\newcommand{\frg}{\mathfrak{g}}
\newcommand{\frh}{\mathfrak{h}}

\newcommand{\frk}{\mathfrak{k}}
\newcommand{\frl}{\mathfrak{l}}
\newcommand{\frn}{\mathfrak{n}}
\newcommand{\fro}{\mathfrak{o}}
\newcommand{\frp}{\mathfrak{p}}
\newcommand{\frq}{\mathfrak{q}}

\newcommand{\frs}{\mathfrak{s}}

\newcommand{\bbC}{\mathbb{C}}

\newcommand{\bbR}{\mathbb{R}}

\newcommand{\bbZ}{\mathbb{Z}}

\newcommand{\caB}{\mathcal{B}}

\newcommand{\caD}{\mathcal{D}}

\newcommand{\caL}{\mathcal{L}}

\newcommand{\caN}{\mathcal{N}}
\newcommand{\caO}{\mathcal{O}}
\newcommand{\caP}{\mathcal{P}}

\newcommand{\caR}{\mathcal{R}}

\newcommand{\caX}{\mathcal{X}}

\begin{document}
\title{regular orbits of symmetric subgroups on partial flag varieties}
\author{Dan Ciubotaru}
\author{Kyo Nishiyama}
\author{Peter E.~Trapa}


\address{Department of Mathematics, University of Utah, Salt Lake City, UT
84112, USA}
\email{ciubo@math.utah.edu}

\address{Department of Mathematics, Graduate School of Science, 
Kyoto University, Sakyo, Kyoto 606-8502, Japan}
\email{kyo@math.kyoto-u.ac.jp}

\address{Department of Mathematics, University of Utah, Salt Lake City, UT
84112, USA}
\email{ptrapa@math.utah.edu}


\maketitle
\definecolor{Red}{rgb}{0,0,0}  
\definecolor{Blue}{rgb}{0,0,0}  
\section{Introduction}

{\color{Red}
The main result of the current paper is a new parametrization
of the orbits of a symmetric subgroup $K$ on a partial flag variety
$\caP$.  The parametrization is in terms of certain Spaltenstein varieties, on one hand, and certain
nilpotent orbits, on the other.   One of our motivations, as explained below,
is related to enumerating special unipotent representations of real reductive groups.
\textcolor{Blue}{Another motivation is understanding (a portion of) the closure order on the set
of nilpotent coadjoint orbits.}

In more detail, suppose $G$ is a complex connected reductive algebraic group
and let $\theta$ denote an involutive automorphism of $G$.
Write $K$ for the fixed points of $\theta$, and $\caP$ for a variety of
parabolic subalgebras of a fixed type in $\frg$, the Lie algebra of $G$.  
Then $K$ acts with finitely
many orbits on $\caP$, and these
orbits may be parametrized in a number of ways (e.g.~\cite{mat}, \cite{rs},
\cite{BH}), each of which may be viewed as a generalization of the classical
Bruhat decomposition.  (This latter decomposition arises
if $G= G_1 \times G_1$, $\theta$
interchanges the two factors, and $\caP$ is taken to be the full flag variety 
of (pairs of) Borel subalgebras.)  We give our parametrization of  $K \bs \caP$
in Corollary \ref{c:param} and then turn to applications and examples in later sections.

\textcolor{Blue}{As mentioned above, one} of the applications we have in mind concerns the connection with nilpotent coadjoint
orbits for $K$.   To each orbit  $Q = K\cdot \frp$ of parabolic subalgebras in $\caP$, we obtain
such a coadjoint orbit as follows.  Let $\frk$ denote the Lie algebra of $K$, and consider
\begin{equation}
\label{e:intro1}
K\cdot\left[(\frg/\frp)^* \cap (\frg/\frk)^*\right ]
\textcolor{Blue}{= K\cdot(\frg/(\frp+ \frk))^* \subset \frg^*};
\end{equation}
here and elsewhere
we implicitly invoke the inclusion of $(\frg/\frp)^*$ and
$(\frg/\frk)^*$ into $\frg^*$ and take the intersection there.
Suppose for simplicity $K$ is connected.  Then the space 
in \eqref{e:intro1} is irreducible.
It also consists of nilpotent elements and is $K$ invariant.
Since the number
of nilpotent $K$ orbits on $(\frg/\frk)^*$ is finite \cite{kr}, the
space must contain a unique dense $K$ orbit, call it
$\Phi_\caP(Q)$.  (It is easy to adapt this argument to yield the
same conclusion if $K$ is disconnected.)
Thus we obtain a natural map
\begin{equation}
\label{e:phi}
\Phi = \Phi_\caP \; : \; K\bs \caP \lra K\bs \caN_\caP^\theta,
\end{equation}
where $\caN_\caP^\theta$ denotes the cone of nilpotent elements
in 
\begin{equation}
\label{e:Ntheta}
\left [G\cdot (\frg/\frp)^*\right ] \cap (\frg/\frk)^*.
\end{equation}
In fact, the map $\Phi_\caP$ is the starting point
of our parametrization of $K\bs \caP$ in Section \ref{s:param}.
For orientation, in the setting of the Bruhat decomposition mentioned above,
the map may be interpretation as taking Weyl group elements to nilpotent
coadjoint orbits.  (Concretely it amounts to taking an element $w$ to the
dense orbit in the \textcolor{Blue}{$G_1$} saturation of the intersection of the nilradicals
of two Borel subalgebras in relative position $w$.)

Just as the Bruhat order on a Weyl group is easier to understand than the classification and
closure order on nilpotent orbits,  the set of $K$ orbits on $\caP$ in some sense  \textcolor{Blue}{behaves
more nicely} than the set of $K$ orbits on $\caN_\caP^\theta$.  \textcolor{Blue}{The former (and the
closure order on it)} can be described
uniformly, for instance \cite{rs}.  This is not the case for $K \bs \caN_\caP^\theta$, where
any (known) classification involves at least some case-by-case analysis.
So a natural question becomes: can one translate the uniform features of
$K$ orbits on $\caP$ to the setting of $K$ orbits on $\caN_\caP^\theta$
using $\Phi_\caP$?  This is the viewpoint we adopt in Section \ref{s:param}.
In particular, one may ask the following: given a $K$ orbit $\caO_K$
in $\caN_\caP^\theta$, does these exist a canonical element $Q$ of $K \bs \caP$
such that $\Phi_\caP(Q) = \caO_K$?  If so, we would be able to embed the set of
$K$ orbits on $\caN_\caP^\theta$ into (the more uniformly behaved) set
of $K$ orbits on $\caP$.  \textcolor{Blue}{One might optimistically hope to
understand a parametrization of $K\bs \caN_\caP^\theta$ (and understand its closure order)
in this way.}

The simplest way to produce  affirmative answer to this last question is if the fiber
of $\Phi_\caP$ over $\caO_K$ consists of a single element $Q$.  So it is desirable
to have a formula for the cardinality of the fiber.  Using ideas of Rossmann and
Borho-MacPherson, we give such a formula in
Proposition \ref{p:count} in terms of certain Springer representations.
The question of whether the fiber consists of a single element then becomes a 
multiplicity one question about certain Weyl group representations.
We then turn to two natural questions: 
\begin{enumerate}
\item[(1)]  Can one find a natural class of 
orbits $\caO_K$
for which the fiber $\Phi_\caP^{-1}(\caO_K)$ is indeed a singleton? 

\item[(2)] If so, can one
give an {\em effective} algorithm to determine the fiber?  (This is clearly
important if one really wants to use these ideas to try to
classify $K$ orbits on $\caN_\caP^\theta$ uniformly.)
\end{enumerate}
We give affirmative answers to these questions in Proposition \ref{p:BiratClosed}
and Remark \ref{r:effective} respectively.
The class of $K$ orbits we find are those $\caO_K$ such that $\caO = G\cdot \caO_K$ is an {\em even}
complex orbit; then 
$\Phi_\caP^{-1}(\caO_K)$ consists of a single element if $\caP$ is taken to be the
partial flag variety  
such that $T^*\caP$ is a resolution
of singularities of the closure of $\caO$.
(The corresponding $K$ orbits 
on $\caP$ are the regular orbits
of the title.)  Perhaps surprisingly
the algorithm answering (2)
relies on the Kazhdan-Lusztig-Vogan algorithm \cite{v:ic3}
for computing the intersection homology groups (with coefficients) 
of $K$ orbit closures on the
full flag variety.

The setting of Section \ref{s:reg} may appear too restrictive to be of
much practical value.  But in Section \ref{s:rt} we recall that it is
exactly the geometric setting of the Adams-Barbasch-Vogan definition
of Arthur packets.    More precisely, since the
ground field is $\bbC$, $\theta$ arises as the complexification of a
Cartan involution for a real form $G_\R$ of $G$.  We show that the
algorithm of Remark \ref{r:effective} gives an effective means to
compute a distinguished constituent of each Arthur packet of integral
special unipotent representations for $G_\R$.  According to the Arthur
conjectures, these representations should be unitary.  This is a
striking prediction (which is still open in general), since the
constructions leading to their definition have nothing to do with
unitarity.

Section \ref{s:rt} is highly technical unfortunately, but we have included it in 
the hope that it is perhaps more
accessible than \cite[Chapter 27]{abv} (upon which it is of course based).   
We have also included it for another reason
which is easy to understand from the current context.  If it were possible to give affirmative
answers to questions (1) and (2) above to a wider class of orbits than we consider in Section \ref{s:reg},
then the ideas of Section \ref{s:rt} translate those answers into new conclusions about
special unipotent representations of real reductive groups.  
In recent joint work with 
Barbasch, one of us (PT) has made progress in this direction. The precise formulation of 
these results involves a rather different set of ideas, and the details will appear elsewhere. 

Finally, in Section \ref{s:example}, we consider a number of examples illustrating some subtleties of the 
parametrization of Section \ref{s:param}. 
}

\noindent{\bf Acknowledgements.} 
KN and PT would like to thank the Hausdorff Research Institute for
Mathematics for 
its hospitality during their stay in 2007.  They would also like to
thank the organizers of the joint MPI-HIM program 
devoted to representation theory, complex analysis and integral geometry.
KN is partially supported by JSPS Grant-in-Aid for Scientific Research (B)
\#{17340037}.
PT was supported by 
NSA grant MSPF-06Y-096 and NSF grant DMS-0532393.

\section{parametrizing $K\bs \caP$}
\label{s:param}

The main result of this section is Corollary \ref{c:param}
which gives a parametrization of the $K$ orbits on $\caP$.
As Propositions \ref{p:count} and \ref{p:AKorbitcount} show, the parametrization
is closely related to Springer's Weyl group representations.

We begin with a discussion of 
the set $K \bs \caB$ of $K$ orbits on $\caB$, \textcolor{Blue}{the full flag variety of 
Borel subalgebras in our fixed complex reductive Lie algebra $\frg$}.
Basic references for this material are \cite{mat} or
\cite{rs}.
The set $K \bs \caB$ is partially
ordered by the inclusion of orbit closures.  It is generated by
closure relations in codimension one.  We will need to distinguish
two kinds of such relations.
To do so, we fix a base-point $\frb_\circ \in \caB$ \textcolor{Blue}{and a Cartan
$\frh_\circ$ in $\frb_\circ$.  We write $\frb_\circ = \frh_\circ \oplus \frn_\circ$
for the corresponding Levi decomposition,
and let $\Delta^+ = \Delta^+(\frh_\circ, \frn_\circ)$ denote the 
roots of $\frh_\circ$ in $\n_\circ$.}
For a simple root $\alpha\textcolor{Blue}{\in \Delta^+}$, let $\caP_\alpha$ denote the set of 
parabolic subalgebras of type $\alpha$, and write $\pi_\alpha$
for the projection $\caB \rightarrow \caP_\alpha$.  

Fix $K$ orbits $Q$ and $Q'$ on $\caB$.  
If $K$ is connected, then $Q$ is irreducible, and hence so is
$\pi_\alpha^{-1}\left (\pi_\alpha(Q)\right)$.  Thus
$\pi_\alpha^{-1}\left (\pi_\alpha(Q)\right)$ contains a unique
dense $K$ orbit.  In general, $K$ need not be connected and
$Q$ need not be irreducible.
But it is easy to see that the similar reasoning applies
to conclude  $\pi_\alpha^{-1}\left (\pi_\alpha(Q)\right)$
always contains a dense $K$ orbit.
We write $Q \stackrel{\alpha}{\rightarrow} Q'$ 
if
\[ \dim(Q') = \dim(Q) +1
\]
 and
\[
\text{$Q'$ is dense in $\pi_\alpha^{-1}\left (\pi_\alpha(Q)\right)$}.
\]
This implies 
that $Q$ is codimension one in the closure
of $Q'$.  The relations $Q < Q'$ for $Q \stackrel{\alpha}{\rightarrow} Q'$ do not
generate the \textcolor{Blue}{full} closure order however.  Instead we must also consider
a kind of saturation condition.  
More precisely, whenever a codimension one
subdiagram of the form
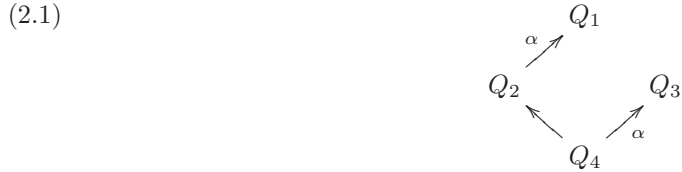
\begin{equation}
\label{e:exchange}
{
\xymatrixcolsep{1pc} 
\xymatrixrowsep{1pc}
\xymatrix
{& Q_1 \\
Q_2 \ar@{->}[ur]^\alpha 
& & Q_3\\
& Q_4 \ar@{->}[ul]\ar@{->}[ur]_\alpha
}
}
\end{equation}
is encountered, we complete it to 
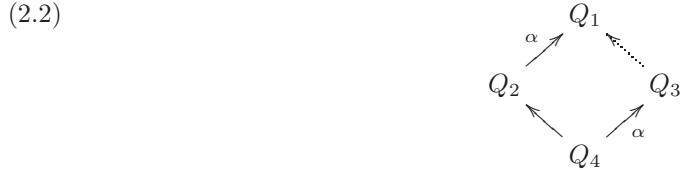
\begin{equation}
\label{e:exchange2}
{
\xymatrixcolsep{1pc} 
\xymatrixrowsep{1pc}
\xymatrix
{& Q_1 \\
Q_2 \ar@{>}[ur]^\alpha 
& & Q_3\ar@{.>}[ul]\\
& Q_4 \ar@{>}[ul]\ar@{>}[ur]_\alpha
}
}
\end{equation}
New edges added in this way are dashed in the diagrams below.  Note
that this operation must be applied recursively, and thus some of the
edges in the original diagram \eqref{e:exchange} may be dashed as the
recursion unfolds.  Following the terminology of
\cite[5.1]{rs}, we call the partially
ordered set determined by the solid edges the weak closure 
order.

Now fix a variety of parabolic subalgebras
$\caP$ of an arbitrary fixed type and write $\pi_\caP$
for the projection from $\caB$ to $\caP$.
For definiteness fix $\frp_\circ \in \caP$
containing $\frb_\circ$, \textcolor{Blue}{and write $\frp_\circ = \frl_\circ \oplus \u_\circ$
for the Levi decomposition such that $\frh_\circ \subset \frl_\circ$}.
Then $K \bs \caP$ may be parametrized from a knowledge of the weak closure
on $K \bs \caB$ as follows.  Consider the relation $Q \sim_\caP Q'$
if $\pi_\caP(Q) = \pi_\caP(Q')$; this is generated
by the relations $Q \sim Q'$ if $Q \stackrel{\alpha}{\rightarrow}Q'$  for
$\alpha$ simple in  $\Delta(\frh_\circ, \frl_\circ)$.
Equivalence classes in $K\bs \caB$
clearly are in bijection with $K\bs \caP$.  
(See also the parametrization of \cite[Section 1]{BH}, especially 
Proposition 4.)
Fix an equivalence class $C$
and fix a representative $Q \in C$.  The same reasoning that
shows that $\pi_\alpha^{-1}\left (\pi_\alpha(Q)\right)$ contains a unique dense
$K$ orbit also shows that
\[
\pi_\caP^{-1}\left ( \pi_\caP(Q) \right )
\]
contains a unique dense $K$ orbit $Q_C \in K\bs \caB$.  In other words, $Q_C$ is the unique
largest dimensional orbit among the elements in $C$.
In fact $Q_C$ is characterized
among the elements of $C$ by the condition
\begin{equation}
\label{e:QC}
\dim \pi_\alpha^{-1}\left (\pi_\alpha(Q_C)\right) = \dim(Q_C)
\end{equation}
for all $\alpha$ simple in $\Delta(\frh_\circ, \frl_\circ)$.
It follows
that the full closure order on $K \bs \caP$ is simply the restriction
of the full closure \textcolor{Blue}{order} on $K\bs \caB$ to the subset of all maximal-dimensional
representatives of the form $Q_C$.  
By restricting only the weak closure order,
we may speak of the weak closure order on $K\bs \caP$.

\textcolor{Red}{
We next place the map $\Phi_\caP$ of \eqref{e:phi} in a more natural context.
\textcolor{Blue}{Consider the cotangent bundle $T^*\caP \subset \caP \times \frg^*$.
It consists of pairs $(\frp, \xi)$ with 
\begin{equation}
\label{e:cotangent}
\xi
\in T_\frp^*\caP \simeq (\frg/\frp)^* 
\end{equation}
The moment map $\mu_\caP$  from $T^*\caP$ to $\frg^*$
maps a point $(\frp, \xi)$ in $T^*\caP$
simply to $\xi$.}  Consider now the
conormal variety for $K$ orbits on $\caP$,
\[
T_K^*\caP = \bigcup_{Q \in K\bs \caP} T^*_Q\caP,
\]
where $T_Q^*\caP$ denotes the conormal bundle to the $K$ orbit $Q$.
(In the special case $G = G_1 \times G_1$ and $\caP = \caB$
mentioned in the introduction,
the conormal variety
is the usual Steinberg variety of triples; see, for instance, the exposition
of \cite{DR}.)
In general we may identify
\begin{equation}
\label{e:conormal}
T_Q^*\caP = \{(\frp, \xi) \; | \; \frp \in Q, \xi \in (\frg/(\frk + \frp))^*\},
\end{equation}
and hence the image of $T_K^*\caP$ under $\mu_\caP$
is simply $\caN_\caP^\theta$.
Moreover, the image of $T_Q^*\caP$ under $\mu_\caP$ is
nothing but the space in \eqref{e:intro1}.  
Hence
$\Phi_\caP(Q)$ is simply
the unique dense $K$ orbit in the moment map image of $T_Q^*\caP$.
}

Here are some
elementary properties of $\Phi_\caP$.

\begin{prop}
\label{p:phi}

\begin{enumerate}
\item
Fix $Q \in K \bs \caP$ and suppose $Q' \in K\bs \caB$ is
dense in $\pi_\caP^{-1}(Q)$.  Then 
\[
\Phi_\caB(Q') = \Phi_\caP(Q).
\]

\item
The map $\Phi_\caP$ is order reversing from the weak closure order in
$K\bs \caP$ to the closure order on $K\bs \caN_\caP^\theta$; that is,
if $Q<Q'$ in the weak closure order on $K \bs \caP$, then
\[
\ol{\Phi_\caP(Q)} \supset \Phi_\caP(Q').
\]
\end{enumerate}
\end{prop}

\noindent {\bf Proof.}  
Part (1)
is clear from the definitions.  Part (2)  reduces to the
assertion for $Q \stackrel{\alpha}{\rightarrow} Q'$.  In that
case, it amounts to a rank one calculation where it is obvious.
\qed

\medskip

\begin{example}
\label{e:sp4full}
Proposition \ref{p:phi}(2) fails for the full closure order on $K\bs \caP$.
The first example which exhibits this failure is $G_\R = \Sp(4,\R)$
and $\caP = \caB$.  Let $\alpha$ denote the short simple root in $\Delta^+$
and $\beta$ the long one.  The closure order for $K \bs \caB$ is
\textcolor{Blue}{as in the diagram in \eqref{e:sp4kgb}}.  
Orbits on the same row of the diagram below all have the
same dimension.  (The bottom row consists of orbits of dimension one,
the next row consists of orbits of dimension two, and so on.)
Dashed lines represent relations in the full closure order which
are not in the weak order.
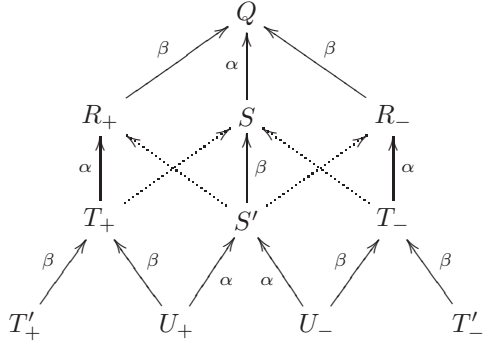
\begin{equation}
\label{e:sp4kgb}
{\xymatrixcolsep{.75pc}
\xymatrixrowsep{2pc}
\xymatrix{
&&&Q\\
& R_+ \ar[urr]^\beta& & S\ar[u]^\alpha  &&
R_-\ar[ull]_\beta \\
& T_+\ar[u]^\alpha \ar@{.>}[urr]_<<<{}& & 
S'\ar[u]_\beta \ar@{.>}[ull]^<<<{} 
\ar@{.>}[urr]_<<<{} &
& T_-\ar[u]_\alpha\ar@{.>}[ull]^<<<{} \\
T'_+ \ar[ur]^\beta & & U_+ \ar[ul]_\beta \ar[ur]_\alpha& &
U_- \ar[ul]^\alpha \ar[ur]^\beta && T'_- \ar[ul]_\beta
}
}
\end{equation}
Adopt the parametrization of $K \bs \caN^\theta$
given in \cite[Theorem 9.3.5]{CM} in terms of signed tableau.
Let $(i_1)^{j_1}_{\epsilon_1}(i_2)^{j_2}_{\epsilon_2}\cdots$ denote the tableau
with $j_k$ rows of length $i_k$ beginning with sign $\epsilon_k$
for each $k$.  Then the closure order on $K\bs \caN^\theta$ is given by
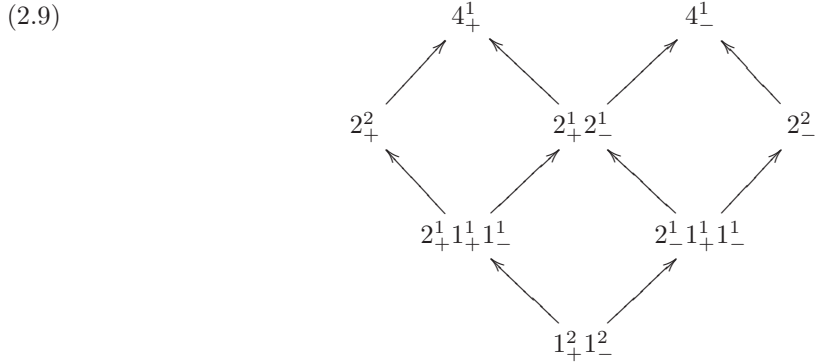
\begin{equation}
\label{e:sp4fullN}
{\xymatrixcolsep{.75pc}
\xymatrixrowsep{2pc}
\xymatrix{
& 4^1_+ \ar@{<-}[dl]\ar@{<-}[dr] && 4^1_- \ar@{<-}[dl]\ar@{<-}[dr]\\
2^2_+\ar@{<-}[dr] & & 2^1_+2^1_-\ar@{<-}[dl]\ar@{<-}[dr]
&& 2^2_-\ar@{<-}[dl]\\
& 2^1_+1^1_+1^1_- \ar@{<-}[dr] & & 2^1_-1^1_+1^1_- \ar@{<-}[dl] \\
&&1^2_+1^2_-\\
}
}
\end{equation}
Then $\Phi_\caB$ maps $Q$ to $1^2_+1^2_-$; $R_\pm$ to $2^1_\pm1^1_+1^1_-$;
$S$ and $S'$ to $2^1_+2^1_-$; $T_\pm$ and $T'_\pm$ to $2^2_\pm$;
and $U_\pm$ to $4^1_\pm$.
Note that $\Phi_\caB$ reverses all closure relations {\em except} the
two dashed edges indicating $T_\pm\subset \ol{S}$.
\end{example}

\medskip

We are now in a position to determine the size of the fiber
$\Phi_\caP^{-1}(\caO_K)$ for $\caO_K \in K\bs \caN_\caP^\theta$.  For
$\xi \in \caO_K$, let $A_K(\xi)$ (resp.~$A_G(\xi)$) 
denote the component group of the centralizer in $K$ (resp.~$G$) 
of $\xi$.   Obviously there is a natural map
\[
A_K(\xi) \rightarrow A_G(\xi)
\]
which we often invoke implicitly.
Write $\Sp(\xi)$
for the Springer representation of $W \times A_G(\xi)$ on
the top homology of the Springer fiber over $\xi$ (normalized
so that $\xi = 0$ gives the sign representations of $W$).  
Let
\[
\Sp(\xi)^{A_K} = \mathrm{Hom}_{A_K(\xi)}\left ( \Sp(\xi), \triv\right).
\]

\begin{prop}
\label{p:count}
Fix $\xi \in \caO_K$.  Then
\begin{align*}
\# \Phi_\caP^{-1}(\caO_K) &= \dim \Hom_{W(\caP)} \left( \sgn, \Sp(\xi)^{A_K} \right )\\
&=\dim \Hom_{W} ( \ind_{W(\caP)}^W(\sgn), \Sp(\xi)^{A_K} ).
\end{align*}
\end{prop}

\noindent{\bf Proof.}  
The second equality follows by Frobenius reciprocity.  For the first, set
\[
S_\caP = \{ Q \in K\bs \caB \; | \; Q \text{ is dense in $\pi_\caP^{-1}
\left( \pi_\caP(Q) \right )$} \}.
\]
According to the discussion around \eqref{e:QC} and Proposition \ref{p:phi}(1),
$\pi_\caP$ implements a bijection
\[
S_\caP \cap \Phi_\caB^{-1}(\caO_K) \rightarrow \Phi_\caP^{-1}(\caO_K).
\]
We will count the left-hand side if $K$ is connected.
If $K$ is disconnected, there are a few subtleties (none of which are
very serious)
which are best treated later.

Consider the top integral Borel-Moore homology
of the conormal variety $T^*_K\caP$.
Since we have assumed $K$ is connected, the closures of the individual conormal bundles exhaust
the irreducible components of $T^*_K\caP$, and their classes form
a basis of the homology,
\[
\H_\top^\infty\left ( T^*_K\caP, \bbZ \right ) = 
\bigoplus_{Q \in K\bs \caP} [\ol{\TQP} ].
\]
If $\caP = \caB$, Rossmann \cite{ross}
(extending earlier work of Kazhdan-Lusztig
\cite{KL})
described a construction giving 
an action \textcolor{Blue}{of the Weyl group} $W$ on this homology space.  The action is graded
in the following sense that if $Q \in \Phi_\caB^{-1}(\caO_K)$, then
\[
w \cdot   [\ol{\TQB} ]
\]
is a linear combination of conormal bundles to orbits in fibers
$\Phi_\caB^{-1}(\caO'_K)$  with ${\caO_K'} \subset \ol{\caO_K}$.  Hence
if we set
\[
\Phi_\caB^{-1}(\caO_K, \leq) = 
\bigcup_{\caO_K' \subseteq \ol{\caO_K}} \Phi_\caB^{-1}(\caO'_K)
\]
and
\[
\Phi_\caB^{-1}(\caO_K, <) = 
\bigcup_{\caO_K' \subsetneq \ol{\caO_K}} \Phi_\caB^{-1}(\caO'_K)
\]
then
\[
\bfM(\caO_K)\; := \; \bigoplus_{Q\in \Phi_\caB^{-1}(\caO_K, \leq)} [\ol{\TQB} ]
\; \;  \biggr / \; 
\bigoplus_{Q\in \Phi_\caB^{-1}(\caO_K, <)} [\ol{\TQB}]
\]
is a $W$ module with basis indexed by $\Phi_\caB^{-1}(\caO_K)$.
Rossmann's construction shows that
\[
\bfM(\caO_K) \simeq \Sp(\xi)^{A_K},
\]
where $\xi \in \caO_K$ as above.  This proves the
proposition for $\caP = \caB$.  For the general case,
we must identify $S_\caP$ in terms of the Weyl group action.
It follows from Rossmann's constructions that
\[
s_\alpha \cdot [\ol{\TQB}] = - [\TQB]
\]
if and only if
\[
\dim \pi_\alpha^{-1}\left (\pi_\alpha(Q)\right) = \dim(Q).
\]
Thus \eqref{e:QC} implies that
$S_\caP \cap \Phi_\caB^{-1}(\caO_K)$ indexes exactly
the basis elements of $\bfM(\caO_K)$ which transform
by the sign representation of the Weyl group of type $\caP$.
The proposition thus follows in the case of $K$ connected.
(A complete proof in the disconnected case is discussed after Proposition 
\ref{p:AKorbitcount}.)
\qed

\medskip

The above proof is extrinsic, in the sense that it is deduced from
a statement about the $\caP = \caB$ case.  We may argue more intrinsically
(without reference to $\caB$) using results of Borho-MacPherson \cite{BM}
as follows.  

Fix $\xi \in \caN_\caP^\theta$ and consider $\muinv_\caP(\xi)$.
In terms of the identification around
\eqref{e:cotangent},
\[
\muinv_\caP(\xi) = \{(\frp,\xi) \; | \; \xi \in (\frg/\frp)^*\}.
\]
(Borho-MacPherson
write $\caP_\xi^0$ for $\muinv_\caP(\xi)$ and call it a Spaltenstein variety.)
 Clearly 
$A_G(\xi)$, and hence $A_K(\xi)$, act on the set of irreducible components
$\Irr(\muinv_\caP(\xi))$.
Fix $C \in \Irr(\muinv_\caP(\xi))$, and consider $Z(C) := \ol{K\cdot C} \subset T^*\caP$.
Since $\xi \in \caN^\theta_\caP \subset \caN(\frg/\frk)^*$, it follows from 
\eqref{e:conormal} that  $Z(C)$
is in fact contained in the conormal variety
\[
Z(C) \subset T^*_K\caP,
\]
which is of course pure-dimensional of dimension $\dim(\caP)$.
Hence
\[
\dim(Z(C)) \leq \dim(\caP).
\]
But clearly
\[
\dim(Z(C))  = \dim(K\cdot\xi) + \dim(C),
\]
and thus
\begin{equation}
\label{e:max}
\dim(C) \leq   \dim(\caP) - \dim(K\cdot \xi).
\end{equation}
Write $\Irr_\max(\muinv_\caP(\xi))$ for those irreducible components whose
dimensions actually achieve the upper bound.  (This set could
be empty, for instance, as we shall see in Example \ref{e:sp4partial}
below when $\caP = \caP_\beta$ and $\xi$ is a representative of
a minimal nilpotent orbit.  \textcolor{Blue}{Note, however, that it is a general theorem of
Spaltenstein's that if $\caP = \caB$, the full flag variety, then
$\Irr_\max(\muinv_\caB(\xi))=\Irr(\muinv_\caB(\xi))$.})

\begin{prop}
\label{p:param}
Fix $\xi \in \caN_\caP^\theta$, set $\caO_K = K\cdot \xi$, 
assume $\Phi_\caP^{-1}(\caO_K)$ is nonempty, and
fix $Q \in \Phi_\caP^{-1}(\caO_K)$.  Then
\[
C(Q) := \ol{\TQP} \cap \muinv_\caP(\xi)
\]
is \textcolor{Red}{the union of elements in an} $A_K(\xi)$ orbit on $\Irr_\max(\muinv_\caP(\xi))$.  The assignment $Q \mapsto C(Q)$
gives a bijection
\begin{equation}
\label{e:parambij}
\Phi_\caP^{-1}(\caO_K) \longrightarrow A_K(\xi) \bs \Irr_\max(\muinv_\caP(\xi)).
\end{equation}
\end{prop}

\noindent {\bf Proof.}
Fix $C \in \Irr_\max(\muinv_\caP(\xi))$.  Then $\dim(Z(C))  = \dim(\caP)$ by definition.
Notice that $Z(C)$ is nearly irreducible (and it is if $K$ is connected).  In general, 
the component group of $K$ (which is
finite by hypothesis) acts transitively on the irreducible
components of $Z(C)$.  But from the definition of $T_K^*\caP$, the closure
of each conormal bundle $\TQP$ consists of a subset of irreducible components
of $T^*_K\caP$ on which the component group of $K$ acts transitively.
Since $\dim(Z(C)) = \dim(T^*_K\caP)$, it follows that there is some 
$Q$ such that
\[
Z(C) = \ol{\TQP};
\]
moreover $Q$ must be an element of $\Phi^{-1}_\caP(\caO_K)$. 
Clearly $Z(C) = Z(C')$ if and only if $C$ and $C'$ are in the same $A_K(\xi)$ orbit.  The assignment $C \mapsto Q$ gives a bijection 
$A_K(\xi) \bs \Irr_\max(\muinv_\caP(\xi)) \rightarrow 
\Phi_\caP^{-1}(\caO_K)$ which, by construction, is the inverse of the
map in \eqref{e:parambij}.
This completes the proof.
\qed

\medskip

\begin{cor}
\label{c:param}
Let $\xi_1,\dots, \xi_k$ be representatives of the $K$ orbits on
$\caN^\theta_\caP$. Then the map
\[
Q \lra \left(\Phi_\caP(Q), \ol{\TQP} \cap \muinv_\caP(\xi_i) \right )
\]
for $i$ the unique index such that $K\cdot \xi_i$ dense in $\Phi_\caP(Q)$
implements a bijection 
\[
K\bs \caP \lra \coprod_i A_K(\xi_i) \bs \Irr_\max(\muinv_\caP(\xi_i)).
\]
\end{cor}

\medskip
Thus everything reduces to understanding the irreducible components of
$\muinv_\caP(\xi)$ of maximal possible dimension.  For this we need 
some nontrivial results of Borho-MacPherson.
\cite[Theorem 3.3]{BM} shows that the fundamental classes
of the elements of $\Irr_\max(\muinv_\caP(\xi))$ index a basis of 
$\Hom_{W(\caP)}(\textcolor{Blue}{\sgn},\Sp(\xi))$.
Actually, to be precise, their condition for $C$ to belong to 
$\Irr_\max(\muinv_\caP(\xi))$ is that
\[
\dim(C) = \dim(\caP) - \frac12 \dim(G\cdot \xi).
\]
To square with \eqref{e:max}, we need to invoke the result of Kostant-Rallis
\cite{kr} that $K\cdot \xi$ is Lagrangian in $G\cdot \xi$.
In any case, because $A_G(\xi)$ acts on $\Sp(\xi)$ and commutes with the $W$ action,
$A_G(\xi)$ also acts on 
$\Hom_{W(\caP)}(\sgn, \Sp(\xi))$, and \cite[Theorem 3.3]{BM}
shows that this action is compatible with the action of $A_G(\xi)$ on 
$\Irr(\muinv_\caP(\xi))$.   In particular this implies the following result.

\begin{prop}
\label{p:AKorbitcount}  
Fix $\xi \in \caN^\theta_\caP$.  Then the number of
 $A_K(\xi)$ orbits on $\Irr_\max(\muinv_\caP(\xi))$ equals the dimension of 
 \[
 \Hom_{W(\caP)} \left( \sgn, \Sp(\xi)^{A_K} \right ).
 \]
\end{prop}

\medskip

\noindent
Combining Proposition \ref{p:param} and \ref{p:AKorbitcount}, we obtain
an alternate proof of Proposition \ref{p:count} which makes no assumption
on the connectedness of $K$.

\begin{remark}
The $\caP = \caB$ case of Corollary \ref{c:param} is due to 
Springer (unpublished).  In this case, $W(\caB)$ is trivial, and
thus $\Phi^{-1}_\caB(\caO_K)$ has order equal to the $W$-representation
$\Sp(\xi)^{A_K}$.
\end{remark}

\medskip

It is of interest to compute the bijection of Corollary \ref{c:param}
as explicitly as possible.  For instance, 
if $G_\R = \GL(n,\C)$ and $\caP = \caB$
consists of pairs of flags, the left-hand side of the 
bijection in Corollary \ref{c:param} consists of elements of
the symmetric group $S_n$.  On the right-hand side, all $A$-groups
are trivial, and the irreducible components in question amount to
pairs of irreducible components of the usual Springer fiber.
Such pairs are parametrized by same-shape pairs of standard
Young tableaux.  Steinberg \cite{stein} showed that the bijection
of the corollary amounts to the classical Robinson-Schensted correspondence.

A few other classical cases have been worked out explicitly
(\cite{vL}, \cite{monty}, \cite{trap:imrn}, \cite{trap:rsbc}).  But
general statements are lacking.  For instance, given $Q$ and
$Q'$, there is no known effective algorithm to decide if $\Phi_\caP(Q)
= \Phi_\caP(Q')$.  The next section is devoted to special cases of the
parametrization which lead to nice general statements.  It might appear
that these special cases are too restrictive to be of much use.  
But it turns out that they encode exactly the geometry needed for the
Adams-Barbasch-Vogan definition of Arthur packets.  This is 
explained in Section \ref{s:rt}.

\section{$\caP$-regular $K$ orbits}
\label{s:reg}

The main results of this section are Proposition \ref{p:BiratClosed}(b)
and Remark \ref{r:effective} which together give an effective
computation of a portion of the bijection of Proposition \ref{p:param}
under the assumption that $\mu_\caP$ is birational.

\begin{defn}[{see \cite[Definition 20.17]{abv}}]
\label{d:large}
\label{d:reg}
A nilpotent orbit $\caO_K$ of $K$ on $\caN_\caP^\theta$
is called $\caP$-regular (or simply regular, if $\caP$ is clear from
the context)
if $G\cdot \caO_K$ is dense in $\mu_\caP(T^*\caP)$.  Since 
$\caO_K$ is Lagrangian in $G\cdot \caO_K$ \cite{kr}, this condition
is equivalent to
\[
\dim(\caO_K) = \frac12 \dim\mu(T^*\caP) = \dim(\frg/\frp),
\]
for any $\frp \in \caP$.
In other
words, $\caP$-regular nilpotent $K$-orbits meet the 
complex Richardson orbit induced from $\frp$.
An orbit $Q$ of $K$ on $\caP$ is called $\caP$-regular 
(or simply regular) if $\Phi_\caP(Q)$
is a $\caP$-regular nilpotent orbit.
Note that regular $\caP$-orbits need not exist in general
(for instance, if $G_\R$ is compact
and $\caP$ is not trivial).
\end{defn}

\medskip

Since regular nilpotent $K$ orbits are automatically maximal in the
closure order on $\caN_\caP^\theta$, Proposition \ref{p:phi}(2) shows
that regular $K$ orbits on $\caP$ are minimal in the weak closure
order:

\begin{prop}
\label{p:weakminimal}
Suppose $Q$ is a regular $K$ orbit on $\caP$.  Then $Q$ is minimal
in the weak closure order on $K \bs \caP$.
\end{prop}

\medskip

The next example shows that regular $K$ orbits on $\caP$ need
not be minimal in the full closure order (i.e.~they need
not be closed).

\begin{example}
\label{e:sp4partial}
Retain the notation of Example \ref{e:sp4full}.  Let 
$\caP_\alpha$ (resp.~$\caP_\beta$) 
consist of parabolic subalgebras of type $\alpha$ (resp.~$\beta$)
and write $\pi_\alpha$ and $\pi_\beta$ in place of $\pi_{\caP_\alpha}$
and $\pi_{\caP_\beta}$, and similarly for $\mu_\alpha$ and $\mu_\beta$.
Then the closure
order on $K\bs \caP_\alpha$ is obtained by
the appropriate restriction from \eqref{e:sp4kgb}.  
(Subscripts now indicate dimensions; dashed edges are those
covering relations present in the full \textcolor{Blue}{closure} order but not the weak one.)
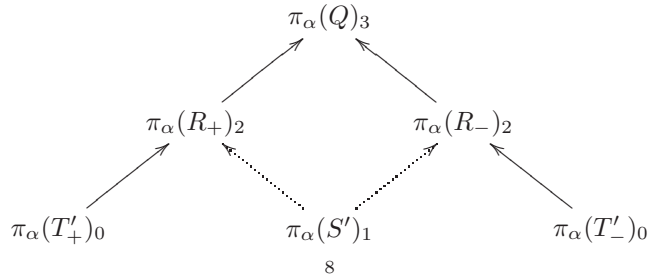
\begin{equation}
\label{e:sp4kgp1}
{\xymatrixcolsep{.75pc}
\xymatrixrowsep{2pc}
\xymatrix{
&&\pi_\alpha(Q)_3\\
& \pi_\alpha(R_+)_2 \ar[ur]
&& 
\pi_\alpha(R_-)_2\ar[ul]
\\
\pi_\alpha(T'_+)_0 \ar[ur] 
&& 
\pi_\alpha(S')_1\ar@{.>}[ul] \ar@{.>}[ur]
&& 
\pi_\alpha(T'_-)_0 \ar[ul]
}
}
\end{equation}
The closure
order on $K\bs \caP_\beta$ is again obtained 
by restriction from \eqref{e:sp4kgb}.  
(Once again subscripts indicate dimensions.)
\begin{equation}
\label{e:sp4kgp2}
{\xymatrixcolsep{.75pc}
\xymatrixrowsep{2pc}
\xymatrix{
&\pi_\beta(Q)_3\\
& \pi_\beta(S)_2 \ar[u]\\
\pi_\beta(T_+)_1 \ar@{.>}[ur] 
&& 
\pi_\beta(T_-)_1 \ar@{.>}[ul] 
}
}
\end{equation}
In this case $\caN^\theta_\alpha = \caN^\theta_\beta$,
and the closure order on $K\bs \caN_\caP^\theta$ is just the bottom
three rows of \eqref{e:sp4fullN},
\begin{equation}
\label{e:sp4partialN}
{\xymatrixcolsep{.75pc}
\xymatrixrowsep{2pc}
\xymatrix
{
2^2_+\ar@{<-}[dr] & & 2^1_+2^1_-\ar@{<-}[dl]\ar@{<-}[dr]
&& 2^2_-\ar@{<-}[dl]\\
& 2^1_+1^1_+1^1_- \ar@{<-}[dr] & & 2^1_-1^1_+1^1_- \ar@{<-}[dl] \\
&&1^2_+1^2_-\\}
}
\end{equation}
From Proposition \ref{p:classical} below (for instance), 
both $\Phi_\alpha = \Phi_{\caP_\alpha}$
and $\Phi_\beta = \Phi_{\caP_\beta}$ are injective.  There are enough
edges in the weak closure order on $K\bs \caP_\alpha$ so that
Proposition \ref{p:phi}(1) allows one to conclude that
$\Phi_\alpha$ reverses the full closure order.  
In fact, $\Phi_\alpha$ is the obvious
order reversing bijection of \eqref{e:sp4kgp1} onto \eqref{e:sp4partialN}.
Hence $\pi_\alpha(T_\pm')$ and $\pi_\alpha(S')$ are $\caP_\alpha$-regular.

By contrast, $\Phi_\beta$ does not invert the dashed edges in 
\eqref{e:sp4kgp2}: $\Phi_\beta$ maps $\pi_\beta(Q)$ to the zero
orbit, and the three remaining orbits to the three orbits of maximal dimension
in $\caN_\caP^\theta$.  Hence $\pi_\beta(T_\pm')$ and $\pi_\beta(S)$ are
$\caP_\beta$-regular.
In particular, $\pi_\beta(S)$ is a $\caP_\beta$-regular
orbit which is not closed.  

Finally note that the fiber of $\Phi_\alpha$
over $2^1_\pm1^1_+1^1_-$ consists of a single element, while the corresponding
fiber for $\Phi_\beta$ is empty.  This is consistent with Proposition 
\ref{p:count} since $\Sp(\xi)$ (for $\xi$ a representative of these orbits)
is a one dimensional representation on which the simple reflection
$s_\alpha$ (resp.~$s_\beta$) acts nontrivially (resp.~trivially).
\qed
\end{example}

\medskip

An essential difference in the two cases considered in Example
\ref{e:sp4partial} is that $\mu_\alpha$ is birational, but $\mu_\beta$
has degree two.

\begin{prop}[{\cite[Theorem 20.18]{abv}}]
\label{p:BiratClosed}
Suppose $\mu_\caP$ is birational onto its image.  Then:
\begin{enumerate}
\item[(a)]
Any regular $K$ orbit on $\caP$
consists of $\theta$-stable parabolic subalgebras (and hence is closed).

\item[(b)]
$\Phi_\caP$ is a bijection from the set of regular $K$ orbits on $\caP$
to the set of regular nilpotent $K$ orbits on 
$\caN_\caP^\theta$.
\end{enumerate}
\end{prop}

\noindent {\bf Proof.}  Fix a $\caP$-regular nilpotent $K$ orbit
$\caO_K$ in $\caN_\caP^\theta$,  $\xi \in \caO_K$,
and $Q \in \Phi_\caP^{-1}(\caO_K)$.  Since $\mu_\caP$ is birational,
the set $\Irr_\max(\mu_\caP^{-1}(\xi))$ is a single point, and so
Proposition \ref{p:param} shows that $Q$ is the unique orbit
in $\Phi^{-1}_\caP(\caO_K)$.  This gives (b).

Again since $\mu_\caP$
is birational, there is a unique parabolic $\frp \in Q$ such that
$\xi \in (\frg/\frp)^*$.  Since $\theta(\xi) =-\xi$, $\theta(\frp)$ is
also such a parabolic.  So $\theta(\frp) = \frp$.  Thus $Q = K\cdot\frp$ 
consists of $\theta$-stable parabolic subalgebras.
This gives the first part of (a).
The same (well-known) proof of the fact that $K$ orbits of $\theta$-stable
Borel subalgebras are closed (for example, \cite[Lemma 5.8]{milicic}), 
also applies to
show that orbits of $\theta$-stable parabolics are closed.  
(It is no longer true that a closed $K$ orbit on $\caP$
consists of $\theta$-stable parabolic subalgebras.  But if a $\theta$-stable
parabolic algebra in $\caP$ exists, 
all closed orbits do indeed consist of $\theta$-stable
parabolic subalgebras.)
\qed 

\medskip

Because of the good properties in Proposition \ref{p:BiratClosed},
we will mostly be interested in $\caP$-regular orbits when 
$\mu_\caP$ is birational.
For orientation (and later use in Section \ref{s:rt})
it is worth recalling a sufficient condition
for birationality from \cite{hesselink}; see also 
\cite[Theorem 7.1.6]{CM} and \cite[Lemma 27.8]{abv}.

\begin{prop}
\label{p:BiratCriterion}
Suppose $\caO$ is an even complex nilpotent orbit.  Let $\caP$ denote
the variety of parabolic subalgebras in $\frg$ corresponding to
the subset of the simple roots labeled 0
in the weighted Dynkin diagram for $\caO$
(e.g.~\cite[Section 3.5]{CM}).  Then $\caO$ is dense
in $\mu_\caP(T^*\caP)$ and $\mu_\caP$
is birational.
\end{prop}
\qed

\medskip

Return to Proposition \ref{p:BiratClosed}(a).  
Example \ref{e:f4} below shows that if $\mu_\caP$ is birational, then
not every (necessarily closed) $K$ orbit of $\theta$-stable parabolic
subalgebras on $\caP$ need be regular.  \textcolor{Blue}{(A good
example to keep in mind is the case when $K$ and $G$ have the same
rank and $\caP = \caB$.  Then the closed $K$ orbits on $\caB$
parametrize discrete series representations with a fixed infinitesimal character.
But the regular orbits are the ones which  parametrize large discrete series.)}
So the question becomes: can
one give an effective procedure to select the regular $K$ orbits on
$\caP$ from among all orbits of $\theta$-stable parabolics (when
$\mu_\caP$ is birational)?  This is only a small part of computing the
parametrization of Corollary \ref{c:param}, so it is perhaps
surprising that the answer we give after Proposition
\ref{p:aql} depends on the power of the
Kazhdan-Lusztig-Vogan algorithm for $G_\R$, the real form
of $G$ with complexified Cartan involution $\theta$.

We need a few definitions.  Recall that
the associated variety of a two-sided ideal $\textcolor{Blue}{I}$
in $\U(\frg)$ is the subvariety of $\frg^*$ cut out by the associated
graded ideal $\gr I$ (with respect to the standard filtration on
$\U(\frg)$) in $\gr\U(\frg) = S(\frg)$.  (From \cite{bb:i},
if $I$ is primitive, then $\AV(I)$ is the closure of  a single nilpotent
coadjoint orbit.)  Finally
if $\frp$ is a $\theta$-stable parabolic subalgebra of $\frg$,
recall the irreducible $(\frg,K)$-module $A_\frp$ constructed in \cite{vz}.
(It would be more customary to denote these modules $A_\frq$, but
we have already used the letter $Q$ for another purpose.)

\begin{prop}
\label{p:aql}
Suppose $\mu_\caP$ is birational.  Fix a closed $K$ orbit $Q$
on $\caP$ consisting of 
$\theta$-stable parabolic subalgebras.  Fix $\frp \in Q$.  Then
$Q$ is $\caP$-regular in the sense of Definition 
\ref{d:large} if and only if
\[
\AV(\Ann(A_\frp)) = \mu(T^*\caP),
\]
the closure of the complex Richardson orbit induced from $\frp$.
\end{prop}

\begin{remark}
\label{r:effective}
We remark that the condition of
the proposition is effectively computable from a knowledge of the
Kazhdan-Lusztig-Vogan polynomials for $G_\R$.  More precisely, the
results of Section \ref{s:param} allow us to enumerate the closed orbits
of $K$ on $\caP$ from the structure of $K$ orbits on
$\caB$.  In turn, the description of $K\bs \caB$ has been implemented
in the command {\tt kgb} in the software package {\tt atlas}
(available for download from {\tt www.liegroups.org}).
\textcolor{Blue}{Moreover,
it is not difficult to determine which closed orbits
consist of $\theta$-stable parabolic subalgebras; in fact, if one of closed orbit
does, then they all do.}
(Alternatively, one may implement the algorithms of \cite[Section 3.3]{BH},
at least if $K$ is connected.)
For a representative
$\frp$ of each such orbit,
one then uses the command {\tt wcells}, to enumerate the cell of
Harish-Chandra modules containing the Vogan-Zuckerman module $A_\frp$.
(The computation of cells relies on computing Kazhdan-Lusztig-Vogan
polynomials.)
Finally $\AV(\Ann(A_\frp)) = \mu(T^*\caP)$ if and only if the cell
containing $A_\frp$ affords the Weyl group representation
$\Sp(\xi)^{A_G}$ (with notation as in Section \ref{s:param}), 
where $\xi$ is an element of
the Richardson orbit induced from $\frp$.
Again, this is an effectively computable condition and is easy to 
implement from the output of {\tt atlas}.  Hence {\em if
  $\mu_\caP$ is birational, there is an effective algorithm to
  enumerate the $\caP$-regular orbits of $K$ on $\caP$.}  
\end{remark}

\begin{remark}
\label{r:effective2}
Suppose $\caO$ is an even complex nilpotent orbit, so that
Proposition \ref{p:BiratCriterion} applies.
Then Proposition \ref{p:BiratClosed}(b) shows that the algorithm
of Remark \ref{r:effective} also
enumerates the $K$ orbits in $\caO \cap (\frg/\frk)^*$.
Using the Kostant-Sekiguchi correspondence, this amounts to the
enumeration of the real forms of $\caO$, i.e.~ $G_\R$ orbits on
$\caO \cap \frg_\bbR^*$.  By contrast, if $\caO$ is not even,
the only known way to enumerate the real forms of $\caO$ involves
case-by-case analysis.
\end{remark}
\medskip

Proposition \ref{p:aql} is known to experts, 
but we sketch a proof (of more
refined results) below; see also \cite[Chapter 20]{abv}.
We begin with some representation-theoretic preliminaries.  Let
$\caD_\caP$ denote the sheaf of algebraic differential operators on
$\caP$, and let $D_\caP$ denote its global section.  Since the
enveloping algebra $\U(\frg)$ acts on $\caP$ by differential
operators, we obtain a map $\U(\frg) \rightarrow D_\caP$.  Let
$I_\caP$ denote its kernel, and $R_\caP$ its image.  By choosing a
base-point $\frp_\circ \in \caP$, it is easy to see that $I_\caP$ is
the annihilator of the irreducible generalized Verma module induced
from $\frp_\circ \in \caP$ with trivial infinitesimal character.
We will be interested in studying Harish-Chandra modules
whose annihilators contain $I_\caP$, i.e.~ $(R_\caP,K)$-modules.  
For orientation, note
that if $\caP = \caB$, $I_\caB$ is a minimal primitive ideal, and
thus any Harish-Chandra module with trivial infinitesimal character
contains it.

Unlike the case of $\caP = \caB$, 
$\U(\frg)$ need not surject onto $D_\caP$ in general,
and so $R_\caP \simeq \U(\frg)/I_\caP$ is generally a proper
subring of $D_\caP$.  Thus the localization functor
\begin{align*}
R_\caP\textrm{-mod} &\lra \caD_\caP\textrm{-mod} \\
X &\lra \caX := \caD_\caP \otimes_{R_\caP} X.
\end{align*}
need not be an equivalence of categories.  
But nonetheless we have that the appropriate irreducible objects match.
(Much more conceptual statements of which the following proposition is a 
consequence have recently been established by S.~Kitchen.)

\begin{prop}
\label{p:res}
Suppose $X$ is an irreducible $(D_\caP, K)$-module.  Then
its restriction to $R_\caP$ is irreducible.
\end{prop}

\noindent{\bf Sketch.}
Irreducible $(D_\caP,K)$-modules are parametrized by
irreducible $K$ equivariant flat connections on $\caP$.  
We show that the irreducible $(R_\caP, K)$-modules
are also parametrized by the same set.  The parametrizations
have the property that support of the localization of either type
of module parametrized by such a  connection $\caL$ is simply the closure
of the support of $\caL$.
This implies there are the same
number of such irreducible modules and hence implies
the proposition.

Let $X$ be an irreducible $(R_\caP, K)$-module.  Hence
we may consider $X$ as an irreducible $(\frg, K)$-module, say $X'$,
whose annihilator contains $I_\caP$.  By localizing on $\caB$,
we may consider the corresponding
irreducible $K$ equivariant flat connection on $\caB$, say $\caL'$, 
parametrizing $X'$. The condition that $\Ann(X') \supset I_\caP$ 
can be translated into a geometric condition on $\caL'$ using
\cite[Lemma 3.5]{lv}, the conclusion of which is that $\caL'$ fibers over
an irreducible flat $K$-equivariant connection on $\caP$ (with fiber
equal to the trivial connection on $\caB_\frl$).  
This implies that irreducible $(R_\caP, K)$-modules
are also parametrized by $K$ equivariant flat connections on $\caP$, as claimed,
and the proposition follows.
\qed
\medskip

\begin{remark}
\label{r:singular}
Proposition \ref{p:res} need not hold when considering  twisted sheaves
of differential operators corresponding to singular infinitesimal characters.
\end{remark}

\medskip 

Next suppose $X$ is an irreducible $R_\caP$ module.  Let
$(\caX^i)$ denote a good filtration on its localization $\caX$
compatible with the degree filtration on $\caD_\caP$.  Let
$\CV(X)$ denote the support of $\gr(\caX)$.  This is well
defined independent of the choice of filtration.  Moreover,
there is a subset $\cv(X) \subset K\bs \caP$ such that
\[
\CV(X) = \bigcup_{Q \in \cv(X)} \ol{\TQP}.
\]
The set $\cv(X)$ is difficult to understand, but there are two
easy facts about it.  First, if $X$ is irreducible, there is a 
dense $K$ orbit, say $\supp(X)$ in the support of $\caX$; then
$\supp(X) \in \cv(X)$.  Moreover if $Q \in \cv(X)$, then
$Q \in \ol{\supp(X)}$.  So, for example, if $\supp(X)$ is closed,
then $\cv(X) = \{\supp(X)\}$.

Finally we define
\[
\AV(X) = \mu(\CV(X)).
\]
(Alternatively one may define $\AV(X)$ as in \cite{v:av} without
localizing.  The fact that the two definitions agree follows
from \cite[Theorem 1.9(c)]{bb:iii}.)  Clearly $\AV(X)$ is the union of closures of $K$
orbits on $\caN_\caP^\theta$.  We let $\av(X)$ denote the set of these
orbits.

Here is how these invariants are tied together.

\begin{theorem}
\label{t:av}
Retain the setting above.  Then
\begin{enumerate}
\item $\AV(I_\caP) = \mu({T^*\caP}).$

\item If $X$ is an irreducible $(R_\caP,K)$-module, then
\[
G\cdot \AV(X) = \AV(\Ann(X)) \subset \AV(I_\caP).
\]
\end{enumerate}
\end{theorem}
\noindent {\bf Proof.}  Part (1) is Theorem 4.6 in \cite{bb:i}.
The equality in part (2) is proved in \cite[Section 6]{v:av};
the inclusion follows because 
$X$ is an $R_\caP = \U(\frg)/I_\caP$ module.
\qed

\begin{prop}
\label{p:cc}
Suppose $X$ is an irreducible $(\caR_\caP, K)$-module
such that there exists a $\caP$-regular $K$ orbit 
$Q\in \cv(X)$.  (For instance, suppose $\supp(X)$
is $\caP$-regular.)
Then $\Phi_\caP(Q)$ is a
$K$ orbit of maximal dimension in $\AV(X)$; 
that is, $\Phi_\caP(Q) \in \av(X)$.
\end{prop}

\noindent {\bf Proof.}  Since $\AV(X) = \mu(\CV(X))$
and since $Q \in \cv(X)$, 
\begin{equation}
\label{e:cc}
\Phi_\caP(Q) \subset \AV(X)
\end{equation}
for any $(R_\caP,K)$-module.  If $Q$ is $\caP$-regular,
then the $G$ saturation of the left-hand side of \eqref{e:cc}
is dense in $\mu(T^*\caP)$.  But by Theorem \ref{t:av}
the right-hand side of \eqref{e:cc} is also contained in 
$\mu(T^*\caP)$.  So the current proposition follows.
\qed

\medskip

\begin{cor}
\label{c:cc}
Suppose $X$ is an irreducible $(R_\caP, K)$-module.
Then the following are equivalent.
\begin{enumerate}
\item[(a)] there exist a $\caP$-regular orbit $Q \in \cv(X)$;

\item[(b)] there exists a $\caP$-regular orbit $\caO_K \in \av(X)$;

\item[(c)] $\Ann(X) = I_\caP$;

\item[(d)] $\AV(\Ann(X)) = \AV(I_\caP)$, 
i.e.~$\AV(\Ann(X)) = \mu(T^*\caP)$.
\end{enumerate}
\end{cor}

\noindent {\bf Proof.}
The equivalence of (a) and (b) follows from the definitions above.
Since the annihilator of any $R_\caP$ module contains $I_\caP$,
the equivalence of (c) and (d) follows from \cite[3.6]{bk}.
Theorem \ref{t:av} and the definitions gives the
equivalence of (b) and (d).
\qed

\medskip

\noindent{\bf Proof of Proposition \ref{p:aql}.}
If $\frp \in \caP$ is a $\theta$-stable parabolic, then the 
Vogan-Zuckerman module $A_\frp$ is the unique irreducible
$(R_\caP, K)$-module whose localization is
supported on the closed orbit 
$K\cdot \frp$ and thus, as remarked
above, $\cv(A_\frp) = \{K \cdot \frp\}$.
So Proposition \ref{p:aql} is a special case of Corollary
\ref{c:cc}.
\qed

\section{applications to special unipotent representations}
\label{s:rt}

The purpose of this section is to explain how the algorithm
of Remark \ref{r:effective} produces special unipotent representations.
Much of this section is implicit in \cite[Chapter 27]{abv}.

Fix a nilpotent adjoint orbit $\caO^\vee$ for $\frg^\vee$, the Langlands
dual of $\frg$.  Fix a Jacobson-Morozov triple $\{e^\vee, h^\vee, f^\vee\}$
for $\caO^\vee$, and set
\[
\chi(\caO^\vee) = (1/2) h^\vee.
\]
Then $\chi(\caO^\vee)$ is an element of some Cartan subalgebra $\frh^\vee$
of $\frg^\vee$.  There is a Cartan subalgebra $\frh$ of $\frg$ such that
$\frh^\vee$ canonically identifies with $\frh^*$.  Hence we may
view
\[
\chi(\caO^\vee) \in \frh^*.
\]
There were many choices made in the definition of $\chi(\caO^\vee)$.  
But nonetheless the infinitesimal character corresponding to
$\chi(\caO^\vee)$ is well-defined; i.e.~$\chi(\caO^\vee)$ is well-defined up
to $G^\vee$ conjugacy and thus (via Harish-Chandra's theorem) specifies
a well-defined maximal ideal $Z(\caO^\vee)$ in the center of $\U(\frg)$.  We call $\chi(\caO^\vee)$
 the unipotent
infinitesimal character attached to $\caO^\vee$.

By a result of Dixmier \cite{dix}, there exists a unique maximal primitive
ideal in $\U(\frg)$ containing $Z(\caO^\vee)$.
Denote it by $I(\caO^\vee)$, and let $d(\caO^\vee)$ denote the dense
nilpotent coadjoint orbit in $\AV(I(\caO^\vee))$.  
The orbit $d(\caO^\vee)$ is called the Spaltenstein dual of $\caO^\vee$
(after Spaltenstein who first defined it in a different way); see
\cite[Appendix A]{bv}.

Fix $G_\R$ as above, and define
\[
\mathrm{Unip}(\caO^\vee) = \{ X \text{ an irreducible $(\frg,K)$ module}
\; | \; \Ann(X) = I(\caO^\vee)\}.
\]
This is the set of special unipotent representations for $G_\R$
attached to $\caO^\vee$.  
Since the annihilator of such a representation $X$ is the maximal primitive
ideal containing $Z(\caO^\vee)$,
$X$ is as small as the (generally singular) infinitesimal character
$\chi(\caO^\vee)$ allows.
These algebraic conditions are conjectured to have
implications about unitarity.

\begin{conj}[Arthur, Barbasch-Vogan {\cite{bv}}]
\label{conj:arthur}
The set $\mathrm{Unip}(\caO^\vee)$ consists of unitary representations.
\end{conj}

\medskip

We are going to produce certain special unipotent representations
from the regular orbits of Definition \ref{d:large}.  
In order to do so, we need to shift our perspective and work on
side of the Langlands dual $\frg^\vee$.  So let $G_\R'$ be a
real form of a connected reductive algebraic group with Lie algebra
$\frg^\vee$ and let $K'$ denote the complexification of a maximal 
compact subgroup in $G'_\R$.  Fix an {\em even} 
nilpotent coadjoint orbit $\caO^\vee$.
(This is equivalent to requiring that $\chi(\caO^\vee)$ is integral.)  
Define $\caP^\vee$ as in Proposition \ref{p:BiratCriterion}.
Thus the main results of Section \ref{s:reg} are available
in this setting.

Let $X'$ denote
an irreducible $(R_{\caP^\vee}, K')$-module,
and let $X$ denote the Vogan dual of $X'$ in the sense of \cite{v:ic4}.  Thus
$X$ is an irreducible Harish-Chandra module for a group $G_\R$ arising
as the real points of a connected reductive algebraic group with Lie
algebra $\frg$.  Moreover, $X$ has trivial infinitesimal character.  

Recall that we are interested in representations with
infinitesimal character $\chi(\caO^\vee)$.  In order to pass to 
this infinitesimal character, we need
to introduce certain translation functors.  
\textcolor{Red}{
There are technical
complications which arise in this setting for two reasons.  First,
$G_\R$ need not be
connected (although it is in Harish-Chandra's class by our hypothesis).
Second, $G_\R$ may not have
enough finite-dimensional representations to define all of the
translations one would like.
Both of these complications disappear if we assume $G$ is simply
connected, and we shall do so here in the interest of streamlining
the exposition.  (It is of course possible
to relax this assumption, as in \cite[Chapter 27]{abv}.)
}

Fix a representative $\rho \in \frh^*$ representing the trivial
infinitesimal character.  Choose a representative $\chi \in \frh^*$ 
representing the (integral) infinitesimal character $\chi(\caO^\vee)$ so that
$\chi$ and $\rho$ lie in the same closed Weyl chamber.  Let 
$\nu = \rho - \chi$.  
Let $F^\nu$ denote the finite-dimensional representation of $G_\R$ with
extremal weight $\nu$; this exists since we have assumed $G$ is simple
connected.  Using it, define the translation functor
$\psi = \psi_\rho^\chi$ (as in \cite[Section VII.13]{kv}) from the category
of Harish-Chandra modules with trivial infinitesimal character to the
category of Harish-Chandra modules with infinitesimal character 
$\chi(\caO^\vee)$.

\begin{theorem}[{cf.~\cite[Chapter 27]{abv}}]
\label{t:sunip}
Retain the notation introduced after Conjecture \ref{conj:arthur}.
In particular, fix an even nilpotent orbit $\caO^\vee$, and let 
$\caP^\vee$ denote the variety of parabolic subalgebras
corresponding to the nodes labeled 0 in the weighted Dynkin
diagram for $\caO^\vee$.
Let $X'$ be an irreducible $(R_{\caP^\vee}, K')$-module,
\textcolor{Red}{
assume $G$ is simply connected,
} 
and 
let $Z = \psi(X)$ denote the translation functor
to infinitesimal character $\chi(\caO^\vee)$ applied to the Vogan dual
$X$ of $X'$.  Then the following are equivalent:
\begin{enumerate}
\item[(a)]  $Z$ is a (nonzero) special unipotent representation
attached to $\caO^\vee$.

\item[(b)] there exists a $\caP^\vee$-regular orbit $Q^\vee \in
\cv(X')$.
\end{enumerate}
\end{theorem}

\noindent{\bf Proof.}
From the properties of the duality explained in \cite[Section 14]{v:ic4}
(and the translation principle),
$Z$ is nonzero with infinitesimal character $\chi(\caO^\vee)$
if and only if $X'$ is annihilated by $I_{\caP^\vee}$, i.e.~if and only
if $X'$ descends to a $(\textcolor{Blue}{R}_{\caP^\vee},K)$-module.  Moreover
$Z$ is
annihilated by a maximal primitive ideal if and only if
the $\textcolor{Blue}{R}_{\caP^\vee}$-module $X'$ has minimal possible annihilator,
namely $I_{\caP^\vee}$.  The conclusion is that 
$Z$ is special unipotent attached to
$\caO^\vee$ if and only if $X'$ is a $(\textcolor{Blue}{R}_{\caP^\vee},K)$-module
annihilated by $I_{\caP^\vee}$.
So the the theorem follows from the equivalence of (a) and (c)
in Corollary \ref{c:cc}.
\qed

\medskip

Since the duality of \cite{v:ic4} is effectively computable,
and since the same is true of the translation functors $\psi$, 
the theorem shows Remark \ref{r:effective} translates into an
effective construction
of special unipotent representations.  More precisely,
one uses Remark \ref{r:effective} to enumerate the relevant 
$\caP^\vee$-regular orbits, and for each one constructs the representation
$X' = A_\frp$ of Proposition \ref{p:aql}.  As remarked in
the proof of Proposition \ref{p:aql}, $X'$ satisfies condition
(b) of Theorem \ref{t:sunip}.  Applying the construction of the theorem
gives special unipotent representations.

In fact, this construction may be understood further in light of the
following refinement.  In the setting of Theorem \ref{t:sunip}, fix a
$\caP^\vee$-regular orbit $Q^\vee$, and define
$\A(Q^\vee)$ be the set of special
unipotent representations attached to $\caO^\vee$ produced by applying
Theorem \ref{t:sunip} to all modules $X'$ with $Q^\vee \in \cv(X')$.
Then the theorem implies
\[
\mathrm{Unip}(\caO^\vee) = \bigcup \A(Q^\vee),
\]
where the (not necessarily disjoint) union is over all $\caP^\vee$-regular
orbits.

The sets $\A(Q^\vee)$ are the Arthur packets defined in \cite[Chapter 27]{abv}.
While there are effective algorithms to
enumerate $\mathrm{Unip}(\caO^\vee)$, there are no such algorithms for individual
packets $\A(Q^\vee)$ (except in favorable cases).
In any event, the discussion of the previous paragraph shows that
{\em Remark \ref{r:effective} leads to an effective algorithm to 
enumerate \textcolor{Blue}{one} element of each Arthur packet of integral special unipotent
representations.}   These representatives are necessarily distinct.

\section{examples}
\label{s:example}

\begin{example}[{\bf Maximal parabolic subalgebras for classical groups}]
\label{e:maximal}
Suppose $G$ is classical and $\caP$ consists of maximal parabolic
subalgebra.  Then it is well-known that 
\[
\ind_{W(\caP)}^W (\sgn)
\]
decomposes multiplicity freely as a $W$-module.  Thus if
$\Sp(\xi)^{A_K}$ is irreducible as a $W$-module, then 
Proposition \ref{p:count} implies $\Phi^{-1}_\caP(\caO_K)$ is a
single orbit.  In particular if the orbits of $A_K(\xi)$ and
$A_G(\xi)$ on irreducible components of the Springer fiber
$\mu_\caB^{-1}(\xi)$ coincide (for instance, if $A_K(\xi)$ surjects
onto $A_G(\xi)$ for each $\xi$), then $\Sp(\xi)^{A_K} =
\Sp(\xi)^{A_G}$ is irreducible and $\Phi_\caP$ is injective.

\begin{prop} 
\label{p:classical}
Suppose the real form $G_\R$ of $G$ corresponding to $\theta$ 
is a classical semisimple Lie group with no
complex factors
whose Lie algebra has no simple factor isomorphic to $\frs\fro^*(2n)$
or $\frs\frp(p,q)$.  If $\caP$ consists of maximal parabolic
subalgebras, then $\Phi_\caP$ is injective.
\end{prop}

\noindent {\bf Proof.}
Unfortunately this follows from a case-by-case analysis
of the classical groups.
First note that the orbits of $A_K(\xi)$ and $A_G(\xi)$
on $\mu_\caB^{-1}(\xi)$ are insensitive to the isogeny class
of $G_\R$.  So, by the remarks preceding the proposition,
it is enough to examine when the two kinds of
orbits coincide for a 
simply connected group $G_\R$ with simple Lie algebra.  
In type A, all $A$-groups are trivial (up to isogeny) so
there is nothing to check.  It follows from direct
computation that $A_K(\xi)$ surjects on $A_G(\xi)$ for
$G_\R = \Sp(2n,\R)$ and $\SO(p,q)$, but that the image
of $A_K(\xi)$ in $A_G(\xi)$ is always trivial for
$\Sp(p,q)$ and $\SO^*(2n)$.  
This completes the case-by-case analysis and hence the
proof.

\medskip

\begin{remark}
For the groups in Proposition \ref{p:classical}, the map $\Phi_\caB$ 
is computed explicitly in \cite{trap:imrn} and \cite{trap:rsbc}.  Using
Proposition \ref{p:phi}(1) this gives one (rather roundabout) way
 to compute $\Phi_\caP$ in these cases.  For exceptional groups,
the injectivity of the proposition fails.  See Example \ref{e:f4}
below.

\end{remark}
\end{example}

\medskip

\begin{example}
\label{e:sp4long}
Suppose now $G_\R = \Sp(2n,\R)$ and $\caP$ consists of maximal
parabolic of type corresponding to the subset of simple roots obtained
by deleting the long one.  (So if $n=2$, 
$\caP = \caP_\alpha$ in Example \ref{e:sp4partial}.)
Then
the analysis of the preceding example extends to show that
 $\Phi_\caP$
is an order-reversing bijection.  The closure order on $K\bs \caN_\caP^\theta$
(and hence $K \bs \caP$) is as follows.
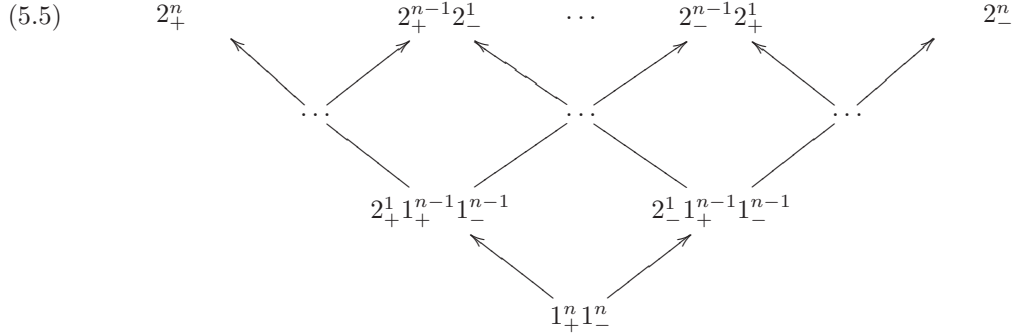
\begin{equation}
\label{e:sp2npartialN}
{\xymatrixcolsep{.75pc}
\xymatrixrowsep{2pc}
\xymatrix
{
2^n_+ \;\;\;\;\;\;\;\;\;\;\ar@{<-}[dr] 
& & 2^{n-1}_+2^1_-\ar@{<-}[dl]\ar@{<-}[dr]
& \dots &
2^{n-1}_-2^1_+\ar@{<-}[dl]\ar@{<-}[dr] &&  \;\;\;\;\;\;\;\;\;\;\;2^n_-\ar@{<-}[dl] \\
&\dots \ar@{-}[dr] & & \dots \ar@{-}[dl]\ar@{-}[dr]
&& \dots\ar@{-}[dl]\\
&& 2^1_+1^{n-1}_+1^{n-1}_- \ar@{<-}[dr] & & 2^1_-1^{n-1}_+1^{n-1}_- \ar@{<-}[dl] \\
&&&1^n_+1^n_-\\}
}
\end{equation}
Here, as before, we are using the parametrization of
$K \bs \caN_\caP^\theta$ given in  \cite[Theorem 9.3.5]{CM}.
There are thus $n+1$ orbits which are $\caP$-regular, all of which
are closed according to Proposition \ref{p:BiratClosed}(a) (which
applies since $\caP$ is attached via Proposition \ref{p:BiratCriterion}
to the even complex orbit with partition $2^n$).

{\color{Blue}
In this setting, we may now apply Theorem \ref{t:sunip}.
(Notationally the roles of the group and dual group must unfortunately
be inverted: for the application, we should take $G^\vee = \Sp(2n,\C)$
in the statement of the theorem.)  
Even though $\SO(n,n+1)$ is not simply connected, the complications
involving the relevant translation functors are absent, and the construction of
the theorem
nonetheless applies and}
produces $n+1$ special unipotent representations for $\SO(n,n+1)$.

\end{example}

\begin{example}
\label{e:unn}
Suppose $G_\R = \U(n,n)$ and $\caP$ corresponds to the subset of simple
roots obtained by deleting the middle simple root in the Dynkin diagram of
type $\mathrm{A}_{2n-1}$.  Then $\Phi_\caP$ is an order reversing bijection,
and the partially ordered sets in question again look like that 
\eqref{e:sp2npartialN} using the parametrization of $K \bs \caN_\caP^\theta$
given in \cite[Theorem 9.3.3]{CM}.
Again there are $n+1$ orbits which are $\caP$-regular.  
\textcolor{Red}{The construction of}
Theorem
\ref{t:sunip} produces $n+1$ special unipotent representation for
$\GL(2n,\R)$, each of which turns out to be a constituent of maximal
Gelfand-Kirillov dimension in the degenerate principal series for $\GL(2n,\R)$
induced from a one-dimensional representation of a Levi factor
isomorphic to a product of $n$ copies of $\GL(2,\R)$.

In terms of representation theory of $G_\R = \U(n,n)$, 
it is well-known that the enveloping algebra in this case
does surject on the ring of global differential operators on $\caP$
(e.g.~the discussion of \cite[Remark 3.3]{trap:aql})
and localization is an equivalence of categories.  Because all Cartan
subgroups in $\U(n,n)$ are connected, the only irreducible flat
$K$-equivariant connections on $\caP$ are the trivial ones supported
on single $K$ orbits.  The map $Q \mapsto \Phi_\caP(Q)$ coincides
with the map which sends the unique irreducible $(R_\caP,
K)$-module supported on the closure of $Q$ to 
the dense orbit in its (irreducible) associated variety, and is
a bijection between such irreducible modules and the $K$
orbits on $\caN_\caP^\theta$.  It would be interesting to see if this
observation could be used to give a geometric explanation of the
computation of composition series of certain degenerate principal series for
$\U(n,n)$ first given in \cite{sahi} and later reproved in \cite{lee}.
(See, for instance, Sahi's module diagrams reproduced in \cite[Figure 7]{lee},
for example.)
\end{example}

\begin{example}
\label{e:sp11}
Suppose $G_\R = \Sp(1,1)$, a real form of $G = \Sp(4,\C)$.
If $\caO$ is the subregular nilpotent orbit for $\frg$ and $\xi \in \caO \cap (\frg/\frk)^*$, then $A_K(\xi)$ is trivial, but $A_G(\xi) \simeq \bbZ/2$.  So the proof
of Proposition
\ref{p:classical} does not apply.  
Let $\alpha$ denote the short simple root and $\beta$ the long one.
The closure order on $K\bs \caB$ is given by
\begin{equation}
\label{e:sp11full}
{
\xymatrixcolsep{1pc} 
\xymatrixrowsep{1pc}
\xymatrix
{& Q\\
& R \ar[u]^{\beta} \\
S_+ \ar[ur]^\alpha&& S_-\ar[ul]_\alpha
}
}
\end{equation}
The picture for $K\bs \caP_\alpha$ 
is
\begin{equation}
\label{e:sp11partialshort}
{
\xymatrixcolsep{1pc} 
\xymatrixrowsep{1pc}
\xymatrix
{\pi_\alpha(Q)_3\\
\pi_\alpha(R)_2 \ar@{->}[u]
}
}
\end{equation}
and for  $K\bs \caP_\beta$
\begin{equation}
\label{e:sp11partiallong}
{
\xymatrixcolsep{1pc} 
\xymatrixrowsep{1pc}
\xymatrix
{& \pi_\beta(Q)_3\\
\pi_\beta(S_+)_2 \ar@{->}[ur]&& \pi_\beta(S_-)_2\ar@{->}[ul]
}
}
\end{equation}
Here $\caN^\theta_\alpha = \caN^\theta_\beta = \caN^\theta_\caB$,
and the closure order of $K$ orbits is simply
\begin{equation}
\label{e:sp11N}
{
\xymatrixcolsep{1pc} 
\xymatrixrowsep{1pc}
\xymatrix
{2^1_+2^1_-\\
1^2_+1^2_- \ar@{->}[u]
}
}
\end{equation}
in the notation of \cite[Theorem 9.3.5]{CM}.  Then $\Phi_\alpha$ is an
order reversing bijection, but $\Phi_\beta$ is two-to-one over
$\textcolor{Blue}{2^1_+2^1_-}$.  The reason is that
\[
\Sp(\xi) = \std \oplus \chi,
\]
where $\std$ is the two-dimensional standard representation of $W$ and
$\chi$ is a character on which the simple reflection $s_\alpha$ acts
trivially and on which $s_\beta$ acts nontrivially.
The orbit $\pi_\alpha(R)$ is $\caP_\alpha$-regular, and the orbits
$\pi_\beta(S_\pm)$ are $\caP_\beta$-regular.
\end{example}

\begin{example}
\label{e:f4}
As an example of what can happen in the exceptional cases, let 
$G$ be the (simply connected) connected complex group of type
$\mathrm{F}_4$ and $\theta$ correspond to the split real form
$G_\R$ of $G$.
(So $K$ is a quotient of
$\Sp(3,\C) \times \SL(2,\C)$ by $\bbZ/2$.)  Then the corresponding
real form $G_\R$ is split.  Let $\caP$ denote the variety of
maximal parabolic obtained by deleting the middle long root from
the Dynkin diagram, and let $\caO$ denote the corresponding
Richardson orbit.  Then $\caO$ is 40 dimensional and
is labeled $\mathrm{F}_4(\mathrm{A}_3)$ in the Bala-Carter
classification.   Moreover $\caO$ is the unique orbit 
which is fixed under Spaltenstein duality.  (Here we are of course identifying $\frg$
and $\frg^\vee$.)
For $\xi \in \caO$,
$A_G(\xi) = S_4$, the symmetric group on four letters.
The weighted
Dynkin diagram of $\caO$ 
has the middle long root labeled 2 and all others
nodes labeled 0.  So $\caP$ corresponds to $\caO$
as in Proposition
\ref{p:BiratCriterion}.

From results of Djokovi\'c (recalled in \cite[Section 9.6]{CM})
there are 19 orbits of $K$ on $\caN_\caP^\theta$.  They are labeled
0--18; the orbit corresponding to label $i$ will be denoted $\caO_K^i$,
and $\xi^i$ will denote an element of $\caO_K^i$.  
Orbits $\caO^{16}_K, \caO^{17}_K$, and  $\caO^{18}_K$ 
are the three $K$ orbits on $\caO \cap 
(\frg/\frk)^*$.  From the discussion leading to \cite[Table 2]{king}, it follows that
$A_K(\xi^i)$ surjects onto $A_G(\xi^i)$ for $i = 0, \dots, 15$.
In each of these cases, $A_G(\xi)$ is either trivial or $\bbZ/2$.
We also have $A_K(\xi^{16}) = A_G(\xi^{16}) = S_4.$  But 
$A_K(\xi^{17}) = D_4$, the dihedral group with eight elements, and
$A_K(\xi^{17}) \rightarrow A_G(\xi^{17})$ is the natural inclusion
into $S_4$.  Finally, $A_K(\xi^{18}) = \bbZ/2 \times \bbZ/2$ which injects
into $A_G(\xi^{18})$.

For $i = 17$ and 18, it is not immediately obvious how to read off
$\Sp(\xi^i)^{A_K(\xi^i)}$ 
from, say, the tables of \cite{carter}.  But
for $i =0, \dots, 16$, the component group calculations of the previous
paragraph imply that $\Sp(\xi^i)^{A_K(\xi^i)} = \Sp(\xi^i)^{A_G(\xi^i)}$, and such
representations are
indeed  tabulated in \cite{carter}.
Applying Proposition \ref{p:count}, it is then not difficult to show
that
\[
\#\Phi^{-1}(\caO_K^i) = 1 \text{ 
if $i \in \{0,1,2,3\} \cup \{9, 10, \dots, 16\}$}
\]
and
\[
\#\Phi^{-1}(\caO_K^i) = 2 \text{ if $i \in \{4, 5,6, 7, 8\}$}.
\]
In more detail, the $G$-saturation of $\caO_K^4$ and $\caO_K^5$
is  the complex orbit $\mathrm{A}_1 \times \wt{\mathrm{A}_1}$ in the Bala-Carter
labeling, while $\caO_K^6$, $\caO_K^7$, and $\caO_K^8$ have $G$ saturation
labeled by $\mathrm{A}_2$.  The corresponding  irreducible Weyl
group representations in these two cases both appear with multiplicity
two in $\ind_{W(\caP)}^W(\sgn)$.  All other relevant multiplicities are one.

We thus conclude that there are 22 orbits of $K$ on $\caP$ which map
via $\Phi_\caP$ to some $\caO_K^i$ for $i = 0, \dots, 15$.
Meanwhile, using the software program {\tt atlas}, one can compute
the closure order of $K$ on $\caB$, and thus (as explained in Section 
\ref{s:param}), the closure order on $K \bs \caP$.  
Figure 4.1
gives the full closure order for $K \bs \caP$.
Vertices are
labeled according to their dimension.  (The edges in Figure
4.1 do {\em not} distinguish between the weak and full closure
order.  Doing so would make the picture significantly more complicated
and difficult to draw.)
There are thus 24 orbits of $K$ on $\caP$.  Since 22 have been shown
to map to $\caO_K^i$ for $i=0,\dots, 15$, one concludes that the
the fiber of $\Phi_\caP$ over $\caO^{i}$ for $i = 16$ and $17$ 
must consist of just one element in each case.

In particular there are three $\caP$-regular $K$ orbits on $\caP$
which are bijectively matched via Proposition \ref{p:BiratClosed}(b)
to $\caO_K^{16}$,
$\caO_K^{17}$, and $\caO_K^{18}$.  But from the {\tt atlas}
computation of the closure order on $K\bs \caP$, there
are {\em four} closed orbits of $K$ on $\caP$.  (These are in
fact exactly the four orbits which are minimal in the weak closure
order.)   See Figure \ref{f:f4}.
The {\tt atlas} labels of the closed orbits are 3, 22, 31, and 47.
Their respective dimensions are 0, 1, 2, and 3.  
Applying the algorithm of Remark \ref{r:effective},
one deduces that the three $\caP$-regular orbits are 3, 31, and 47.  
Theorem \ref{t:sunip}
thus produces three distinct special unipotent representations,
one in each of the three Arthur packets for $\caO = d(\caO)$.

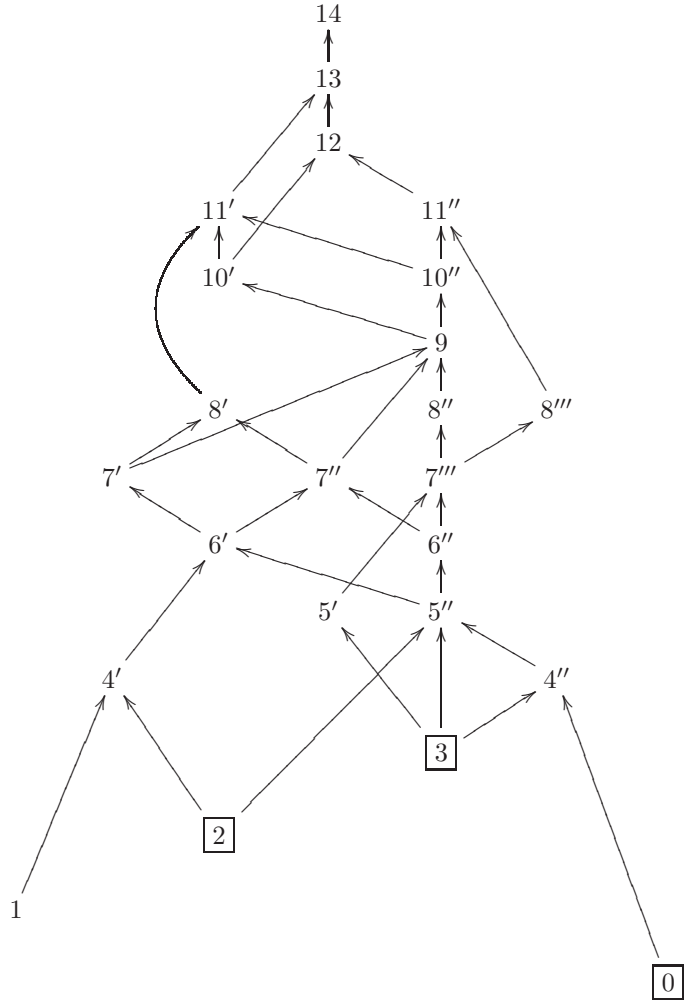
\begin{figure}
\label{f:f4}
\begin{equation*}
{
\xymatrixcolsep{2pc} 
\xymatrixrowsep{1pc}
\xymatrix
{
&&&&14 \\
&&&&13\ar@{->}[u]\\
&&&&12\ar@{->}[u]\\
&&&11'\ar@{->}[uur]&&11''\ar@{->}[ul]\\
&&&10'\ar@{->}[u]\ar@{->}[uur]
&&10''\ar@{->}[u]\ar@{->}[ull]\\
&&&&&9\ar@{->}[u]\ar@{->}[ull]\\
&&&8'\ar@/^{2pc}/[uuu]&&8''\ar@{->}[u] & 8'''\ar@{->}[uuul]\\
&&7' \ar[ur] \ar[uurrr] & & 7''\ar[ul]\ar[uur] & 7'''\ar[u]\ar[ur]\\
&&&6' \ar[ul]\ar[ur] && 6'' \ar[ul]\ar[u]\\
&&&&5' \ar[uur]& 5''\ar[ull]\ar[u] \\
&& 4'\ar[uur]&&&&4''\ar[ul]\\
&&&&&\boxed{3}\ar[uul]\ar[uu]\ar[ur]\\
&&&\boxed{2}\ar[uul]\ar[uuurr]\\
&1\ar[uuur]\\
&&&&&&&\boxed{0}\ar[uuuul]
}
}
\end{equation*}
\caption{The full closure ordering of $K$-orbits on $\caP$ for 
$G_{\R}=\mathrm{F}_4$ and $\caO =\mathrm{F}_4(\mathrm{A}_3)$.  Vertices
are labeled according to their dimensions and boxed vertices
are $\caP$-regular.  Note, in particular, that not every closed
orbit is $\caP$-regular.}
\end{figure}

\end{example}

\end{document}